\documentclass[sn-mathphys]{sn-jnl}
\usepackage{amsmath}
\usepackage{amssymb}
\usepackage{amsfonts}
\usepackage{amsthm}
\usepackage{bbm}
\usepackage{wrapfig}
\usepackage{bm}
\usepackage[mathscr]{eucal}
\usepackage{geometry}
\usepackage{graphicx}
\usepackage{color}
\usepackage{hyperref}
\usepackage{chngcntr}

\counterwithout{figure}{section}

\theoremstyle{thmstyleone}
\newtheorem{theorem}{Theorem}
\newtheorem{proposition}{Proposition}

\newtheorem{innercustomthm}{Theorem}
\newenvironment{customthm}[1]
  {\renewcommand\theinnercustomthm{#1}\innercustomthm}
  {\endinnercustomthm}

\newtheorem{innercustomprop}{Proposition}
\newenvironment{customprop}[1]
  {\renewcommand\theinnercustomprop{#1}\innercustomprop}
  {\endinnercustomprop}

\begin{document}

\title{Diameter of compact Riemann surfaces}

\author*[1,2,4]{\fnm{Huck} \sur{Stepanyants,}}\email{stepanyants.j@northeastern.edu}

\author[3]{\fnm{Alan} \sur{Beardon,}}\email{afb@dpmms.cam.ac.uk}

\author[1,2]{\fnm{Jeremy} \sur{Paton,}}\email{paton.j@northeastern.edu}

\author*[1,2,4,5]{\fnm{Dmitri} \sur{Krioukov}}\email{dima@northeastern.edu}

\affil[1]{\orgdiv{Department of Physics}, \orgname{Northeastern University}, \orgaddress{\city{Boston}, \postcode{02115}, \state{Massachusetts}, \country{USA}}}

\affil[2]{\orgdiv{Network Science Institute}, \orgname{Northeastern University}, \orgaddress{\city{Boston}, \postcode{02115}, \state{Massachusetts}, \country{USA}}}

\affil[3]{\orgdiv{Department of Pure Mathematics and Mathematical Statistics}, \orgname{University of Cambridge}, \orgaddress{\city{Cambridge}, \postcode{CB3 0WB}, \country{UK}}}

\affil[4]{\orgdiv{Department of Mathematics}, \orgname{Northeastern University}, \orgaddress{\city{Boston}, \postcode{02115}, \state{Massachusetts}, \country{USA}}}

\affil[5]{\orgdiv{Department of Electrical and Computer Engineering}, \orgname{Northeastern University}, \orgaddress{\city{Boston}, \postcode{02115}, \state{Massachusetts}, \country{USA}}}

\abstract{Diameter is one of the most basic properties of a geometric object, while Riemann surfaces 
are one of the most basic geometric objects. Surprisingly, the diameter of compact Riemann surfaces 
is known exactly only for the sphere and the torus. For higher genuses, only very general but 
loose upper and lower bounds are available. The problem of calculating the diameter 
exactly has been intractable since there is no simple expression for the distance between a 
pair of points on a high-genus surface. Here we prove that the diameters of a class of simple Riemann 
surfaces known as \emph{generalized Bolza surfaces} of any genus greater than $1$ are 
equal to the radii of their fundamental polygons. This is the first exact result for the 
diameter of a compact hyperbolic manifold.}


\maketitle

\section{Introduction}\label{sec1}

A Riemann surface is any connected, one-dimensional complex manifold. According to the 
classification theorem of closed surfaces~\cite{Seifert}, 
any compact Riemann surface is homeomorphic to either the sphere or the connected sum of $g$ tori, where 
$g$ is the genus of the surface. Once equipped with a metric, a Riemann surface becomes a 
Riemannian manifold. The sphere $(g = 0)$ admits the spherical metric, the torus $(g = 1)$ admits the 
Euclidean metric, while the surfaces of genus $g > 1$ admit the hyperbolic metric. 

The diameter $\mathscr{D}$ of a metric space is the maximum distance between a pair of points in it. 
Diameter is one of the 
most basic characteristics of any geometric object. Surprisingly, the diameter of compact Riemann surfaces is 
known exactly only for the sphere and the torus. The best results on the 
diameter of surfaces of genus $g > 1$ are only loose lower and upper bounds in terms of the total surface 
area, the systole, and the genus~\cite{Chang,Balacheff,Budzinski,Bavard}. 
The \emph{systole} is the length of a shortest noncontractible loop on the surface. 

Our main result is a proof of the following theorem: 
\begin{theorem}\label{T1}
    The diameter of a class of Riemann surfaces~$S_g$, known as generalized Bolza surfaces, of genus 
    $g > 1$ is $\mathscr{D}_g = \textnormal{arccosh} \left( \cot^2 (\pi /4g) \right)$. 
\end{theorem}

The definition of the surfaces $S_g$ is in the next section. In particular, $S_2$ is 
known as the Bolza surface~\cite{Bolza}, one of the first compact hyperbolic manifolds ever considered. 
A free particle moving along a geodesic on the 
Bolza surface was the first dynamical system proven rigorously to be chaotic~\cite{Hadamard}. 
The Bolza surface is 
also known to maximize the systole across all genus-$2$ surfaces~\cite{Werner}. 
For ${g>2}$, the $S_g$ are the generalized Bolza 
surfaces~\cite{Ebbens}. They appear frequently in studies of hyperbolic surfaces due to their 
high degree of symmetry~\cite{Wiman,Ebbens,Bujalance1,Bujalance2}.

Knowing the diameter of a surface, we can efficiently compute the distance between any pair of points 
on it, a result to be published elsewhere. Theorem~\ref{T1} 
states that the Bolza surface has diameter 
$\mathscr{D} = \textnormal{arccosh} \left( 3 + 2 \sqrt{2} \right) \approx 2.45$. 
To the best of our knowledge, this is the first ever exact result for the diameter of \emph{a} 
compact hyperbolic manifold. The closest results to ours appear to be the ones in~\cite{Nabutovsky2003}.
They apply to manifolds of dimension at least five.

We proceed by collecting all the necessary background information and definitions in Section \ref{S2}. 
Section \ref{S3} contains the outline of the proof of Theorem~\ref{T1} split into a sequence of theorems that 
we state in that section as well. We prove all those theorems in the concluding Section \ref{S4}. 

\section{Background information and definitions} \label{S2}

\noindent We use the Poincar\'e disk model of the hyperbolic plane $\mathbb{H}^2$. 
The isometries (distance-preserving maps) of 
$\mathbb{H}^2$ are given by the matrices 
\begin{equation} \label{E1}
    \left[ \begin{matrix}
        a & \overline{c} \\ c & \overline{a}
    \end{matrix} \right], \quad
    a,c \in \mathbb{C}, ~
    \lvert a \rvert ^2 - \lvert c \rvert ^2 = 1,
\end{equation}
which form a subgroup of $PSL(2, \mathbb{C})$, the invertible 
${2 \times 2}$ complex matrices. The action of each matrix on a complex number 
$z$ in the Poincar\'e disk~${\mathbb{D} = \{ z \in \mathbb{C} : \lvert z \rvert < 1\}}$ 
is a fractional linear transformation: 
\begin{equation}
    \left[ \begin{matrix}
        a & \overline{c} \\ c & \overline{a}
    \end{matrix} \right]
    (z) =
    \frac{az + \overline{c}}{cz + \overline{a}}.
\end{equation}

A \emph{Fuchsian group} $\mathscr{F}$ is a discrete subgroup of ${PSL(2,\mathbb{C})}$ 
that has an invariant disk in ${\mathbb{C}^\infty}$. Each $\mathscr{F}$ 
defines a hyperbolic Riemann surface $S$ which is the quotient surface 
${ S = \mathbb{H}^2/\mathscr{F} }$. A 
\emph{fundamental domain} is an open, connected set in~${\mathbb{H}^2}$ that 
contains at most one representative of each point on ${S}$, and whose closure contains
at least one representative of each point on ${S}$~\cite{Beardon}. A 
fundamental domain and its images under the actions of~$\mathscr{F}$ tessellate~$\mathbb{H}^2$.

Since points on a quotient surface $S$ are cosets in $\mathbb{H}^2/\mathscr{F}$, the distance 
between two points (cosets) $[z] = \mathscr{F} z, [w] = \mathscr{F} w$ on $S$ is given by 
\begin{equation} \label{Edstar}
    \delta^\star ([z], [w]) = \inf \{ \delta(z^\prime, w^\prime), z^\prime \in [z], w^\prime \in [w] \} ,
\end{equation}
where $\delta (z, w)$ denotes the distance in $\mathbb{H}^2$. 
The diameter $\mathscr{D}$ of $S$ is the largest distance 
between two points on $S$. 

The best results on the diameters of 
Riemann surfaces of genus $g > 1$ are as follows. It was shown in~\cite{Chang} that for any such 
surface, the following inequalities hold, where $\ell$ 
is the systole, $A$ is the total area of the surface, and $\mathscr{D}$ is the diameter: 
\begin{align}
    2 \ell \sinh (\mathscr{D}) & \geq A \text{, and } \\
    2 \sinh(\ell /4) \mathscr{D} & \leq A .
\end{align}
In~\cite{Balacheff} it was shown that 
\begin{equation}
    4 \cosh (\ell /2) \leq 3 \cosh (\mathscr{D}) - 1 .
\end{equation}
Another lower bound exists in terms of the area alone~\cite{Budzinski}:
\begin{equation}
    \cosh \mathscr{D} \geq \frac{A}{2 \pi} + 1 ,
\end{equation}
and finally, there is the following lower bound in terms of the genus~\cite{Bavard}:
\begin{equation}
    \cosh \mathscr{D} \geq \frac{1}{\sqrt{3}} \cot \left( \frac{\pi}{6(2g-1)} \right) .
\end{equation}

All of the bounds above hold in general for any Riemann surface, as do other 
related spectral results~\cite{Cheng1,Cheng2,Cheeger}. However, none of them is tight 
for the~${S_g}$.

\begin{figure}[t!]
    \begin{center}
        \includegraphics[width=0.35\textwidth]{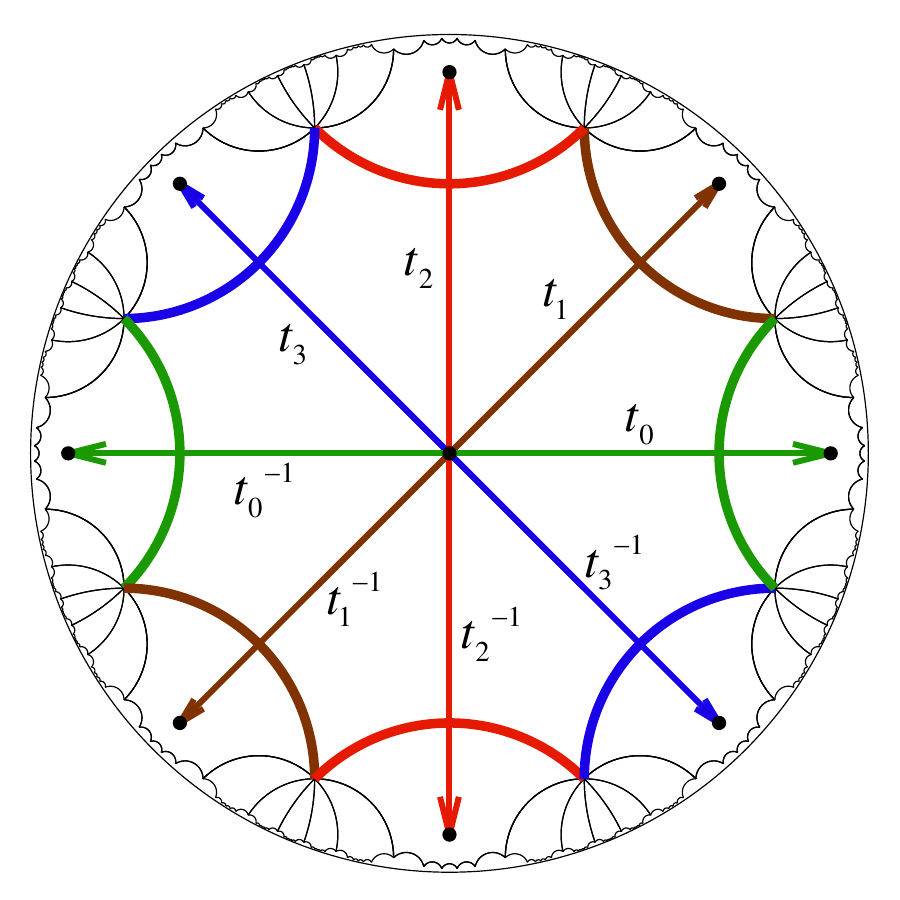}
        \includegraphics[width=0.35\textwidth]{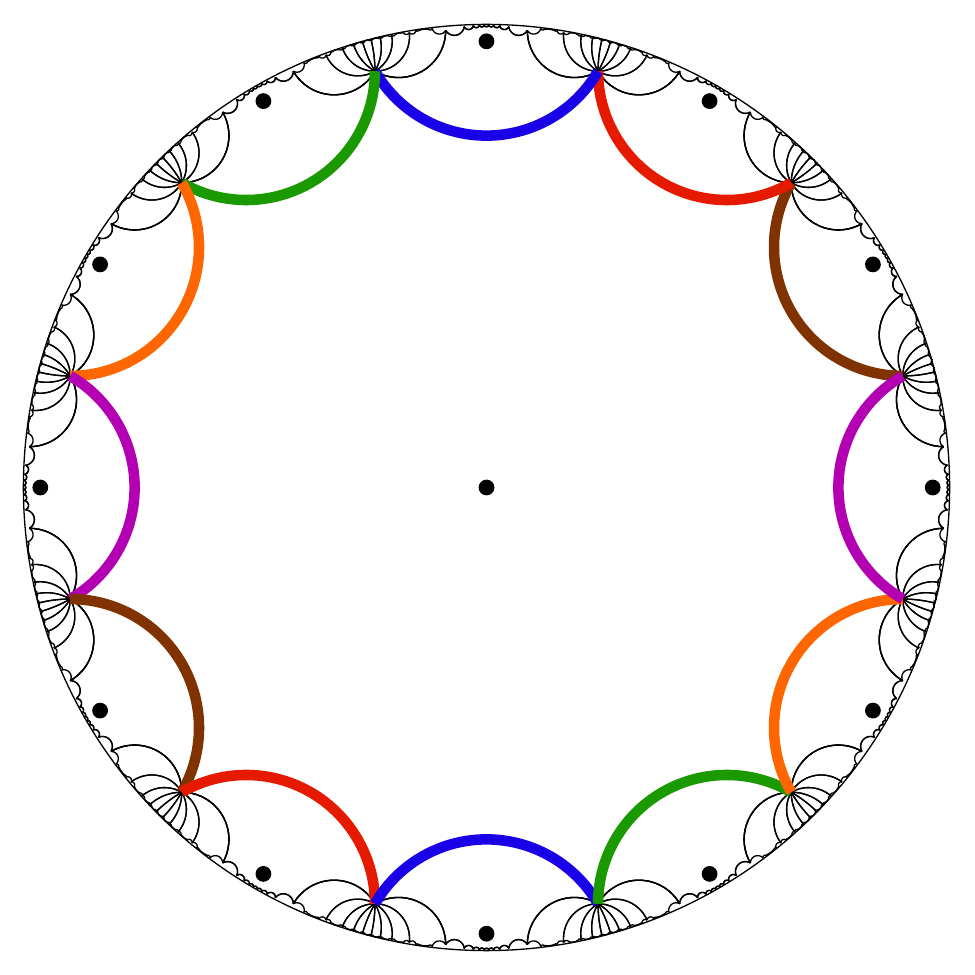}
        \caption{Left: The Bolza surface $(S_2)$ and its ``gluing'' scheme. Sides that 
                 are identified in the quotient space are labeled with the same color, 
                 and the arrows represent the actions of the generators ${ \{t_k\} }$ 
                 and their inverses on the fundamental polygon. 
                 Right: $S_3$ and its gluing scheme.}
        \label{F1}
    \end{center}
\end{figure}

The \emph{generalized Bolza surface} ${S_g}$ of genus ${g}$ is 
defined~\cite{Ebbens,Aurich} to be the surface obtained by identifying
(``gluing'') the opposite sides of the regular $4g$-gon with interior angles $\pi/2g$,
Fig.~\ref{F1}, whose side length~$s$ and radius~$R$ satisfy~\cite{Aurich}
\begin{equation} \label{EsR}
    \cosh ^2 \left( \frac{s}{2} \right) = \cosh R = \cot^2 \left( \frac{\pi}{4g} \right),
\end{equation}
and whose vertices $\{ v_k \}$ are evenly spaced at distance $R$ from the origin: 
\begin{equation}
    v_k = \tanh \left( \frac{R}{2} \right) ~ e^{\left( k - \frac{1}{2} \right) \frac{\pi i}{2g}} ,
    \quad k = 0,1, \dots, 4g - 1 .
\end{equation}
This surface~$S_g$ is the quotient surface ${\mathbb{H}^2 / \mathscr{F}_g}$,
where ${\mathscr{F}_g}$ is the Fuchsian group generated by the following ${2g}$
generators and their inverses:
\begin{equation}
    t_k = 
    \left[ \begin{matrix} 1 & \tanh \left( \frac{s}{2} \right) e^{\frac{k \pi i}{2g}} \\ 
    \tanh \left( \frac{s}{2} \right) e^{- \frac{k \pi i}{2g}} & 1\end{matrix} \right], 
    \quad k = 0,1, \dots, 2g - 1.
\end{equation}
These generators glue the opposite sides of the polygon by mapping it to its edge-adjacent
polygons in the tessellation as shown in Fig.~\ref{F1}. 
By Poincar\'e's theorem~\cite{Beardon}, the interior of the polygon is a fundamental domain
of~${S_g}$, while the polygon itself is called the \emph{fundamental polygon}.

The surface~$S_2$ of genus $2$ is the well-known 
Bolza surface~\cite{Bolza}. The surfaces~${S_g}$ are also known as 
the \emph{Wiman surfaces of type II}~\cite{Bujalance1,Bujalance2}. 
They have ${8g}$ automorphisms, except 
the ${g=2}$ Bolza surface, which has ${48}$ automorphisms. 
The ${S_g}$ also have automorphisms of order ${4g}$, which 
is the second largest possible order, as was proven in~\cite{Wiman}.

Our main result is the exact diameter of the surface~${S_g}$ of genus ${g > 1}$, 
Theorem~\ref{T1}, which in view of Eq.~(\ref{EsR}) can be restated as: 
\begin{customthm}{1} \label{T1p}
    ${\mathscr{D}_g = R}$. 
\end{customthm}

The problem of finding the diameter 
is not trivial thanks to the definition of the distance in Eq.~(\ref{Edstar}), 
which is an infimum over the infinitely many elements in~${\mathscr{F}_g}$.
Even though this infimum is actually a minimum since only a finite subset of
group elements needs to be considered, there is no explicit expression for
such a subset, which is the main difficulty in the problem. In fact, as mentioned in the introduction,
our main motivation for this paper is that the knowledge of the diameter of a surface
allows us to provide an explicit expression for this subset, leading to an efficient formula
to compute distances between pairs of points on the surface, a result to appear in a follow-on paper.

Our proof of the main result, which we outline in the following section, is a 
combination of geometric and algebraic techniques. 
Specifically, we use some geometric symmetries of $S_g$ to simplify the problem, 
and algebra to compute the optimal distances. 

\section{Proof strategy}\label{S3}

\noindent We first observe that 
every vertex of the fundamental polygon represents the same point in the quotient 
space; we will call this point the \emph{quotient vertex}, $[v] \in S_g$. 
Similarly, we define $[0]$ as the point in the quotient space to which the origin maps. 
For an arbitrary point $z$ in the Poincar\'e disk, let $d_0(z)$ and $d_v(z)$ denote 
the quotient distances from $[z]$ to $[0]$ and $[v]$, respectively. 

The proof of Theorem~\ref{T1} will be broken into several smaller, sequential theorems. 
First, we will show that ${\delta^\star([0],[v]) = R}$, thus proving the following 
theorem: 

\begin{theorem} \label{T2}
    The diameter $\mathscr{D}_g$ of $S_g$ is at least $R$, 
    $\mathscr{D}_g \geq R$. 
\end{theorem}

The remaining theorems are aimed at proving $\mathscr{D}_g \leq R$, 
or equivalently, ${\delta^\star ([z],[w]) \leq R}$ for 
every pair of points ${[z],[w]}$ in the quotient space. 
The symmetry of $S_g$ allows us to make several simplifying 
assumptions about $[z]$ and $[w]$ without loss of generality. 

The symmetries of the quotient surface $S_g$ are 
isometric bijections that map ${S_g}$ to itself. In the Poincar\'e disk ${\mathbb{D}}$, 
these symmetries can be 
thought of as the subgroup of isometries $\phi$ of ${\mathbb{D}}$ for which 
the projection ${\phi^\star}$ of ${\phi}$ onto the quotient surface, given by
\begin{equation}
    \phi^\star([p]) = [\phi(p)] ,
\end{equation}
is well-defined and isometric. These properties are satisfied by all ${\phi}$ 
for which the Fuchsian group ${\mathscr{F}_g}$, a subgroup of ${PSL(2,\mathbb{C})}$, is invariant under 
conjugation by ${\phi}$: ${\phi \mathscr{F}_g \phi^{-1} = \mathscr{F}_g}$. Indeed, for such 
${\phi}$ and any ${f,f_1,f_2 \in \mathscr{F}_g}$ and ${p,p_1,p_2 \in \mathbb{D}}$, there exist 
${f^\prime,f_1^\prime,f_2^\prime \in \mathscr{F}_g}$ such that
\begin{align}
    & \phi^\star([fp]) = [\phi f(p)] = [f^\prime \phi (p)] = [\phi(p)] = \phi^\star([p]), \text{ and} \\
    \nonumber
    & \delta^\star(\phi^\star([p_1]), \phi^\star([p_2])) = 
    \delta^\star([\phi (p_1)], [\phi (p_2)]) = \\
    \nonumber
    & \min \{ \delta(f_1 \phi(p_1),f_2 \phi(p_2)), f_1,f_2 \in \mathscr{F}_g \} = \\
    \nonumber
    & \min \{ \delta(\phi f_1^\prime (p_1),\phi f_2^\prime (p_2)), f_1^\prime,f_2^\prime \in \mathscr{F}_g \} = \\
    \nonumber
    & \min \{ \delta(f_1^\prime (p_1)), \delta(f_2^\prime (p_2)), f_1^\prime,f_2^\prime \in \mathscr{F}_g \} = \\
    & \delta^\star([p_1],[p_2]),
\end{align}
so ${\phi}$ is well-defined and an isometry. 

The simplest such isometries are the elements of $\mathscr{F}_g$. Indeed, for 
${\phi \in \mathscr{F}_g}$ we have ${\phi \mathscr{F}_g \phi^{-1} = \mathscr{F}_g}$. 
In this case, ${\phi^\star}$ is simply the identity function on ${S_g}$. 
Therefore, $S_g$ is symmetric under actions by elements of ${\mathscr{F}_g}$.

Another subgroup of isometries of $\mathbb{D}$ that are symmetries of $S_g$ are 
rotations of $\mathbb{D}$ by ${\pi / 4g}$. Although 
this fact is evident from the definition of 
$S_g$ in Section~\ref{S2}, we give its formal proof in Section~\ref{S4}. 

\begin{figure}[t!]
    \centering
    \includegraphics[width=0.3\textwidth]{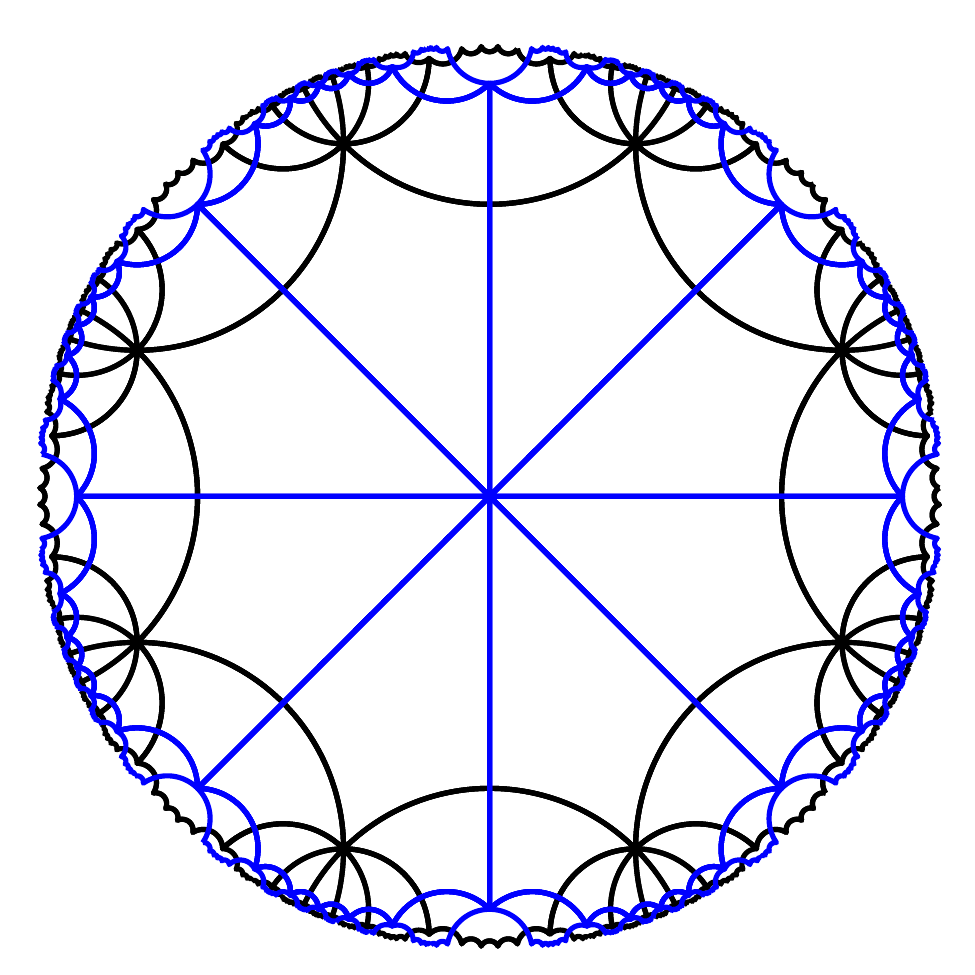}
    \caption{The dual tessellation (blue) overlayed against the original tessellation (black) in the 
             case $g = 2$.}
    \label{F2}
\end{figure}

The final and most complicated subgroup of symmetries of $S_g$ 
that we will need is its \emph{dual symmetries}. 
We define them with the aid of new fundamental polygons that 
we call the \emph{dual polygons}, 
and the corresponding \emph{dual tessellation} of $\mathbb{H}^2$. Dual polygons are 
formed by joining by geodesics the centers of the $4g$ polygons around the 
vertices in the original tessellation, Fig.~\ref{F2}. We will prove 
the following proposition: 

\begin{proposition} \label{P1}
    Let ${\phi}$ be any isometry of $\mathbb{D}$ that maps the fundamental polygon to 
    any of its dual polygons. Then 
    the projection $\phi^\star : S_g \rightarrow S_g$ given by 
    \begin{equation}
        \phi^\star ([z]) = [\phi(z)]
    \end{equation}
    is well-defined and is an isometric automorphism of $S_g$. In particular, 
    dual polygons are fundamental polygons. 
\end{proposition}

\begin{figure}[t!]
    \begin{center}
        \includegraphics[width=0.35\textwidth]{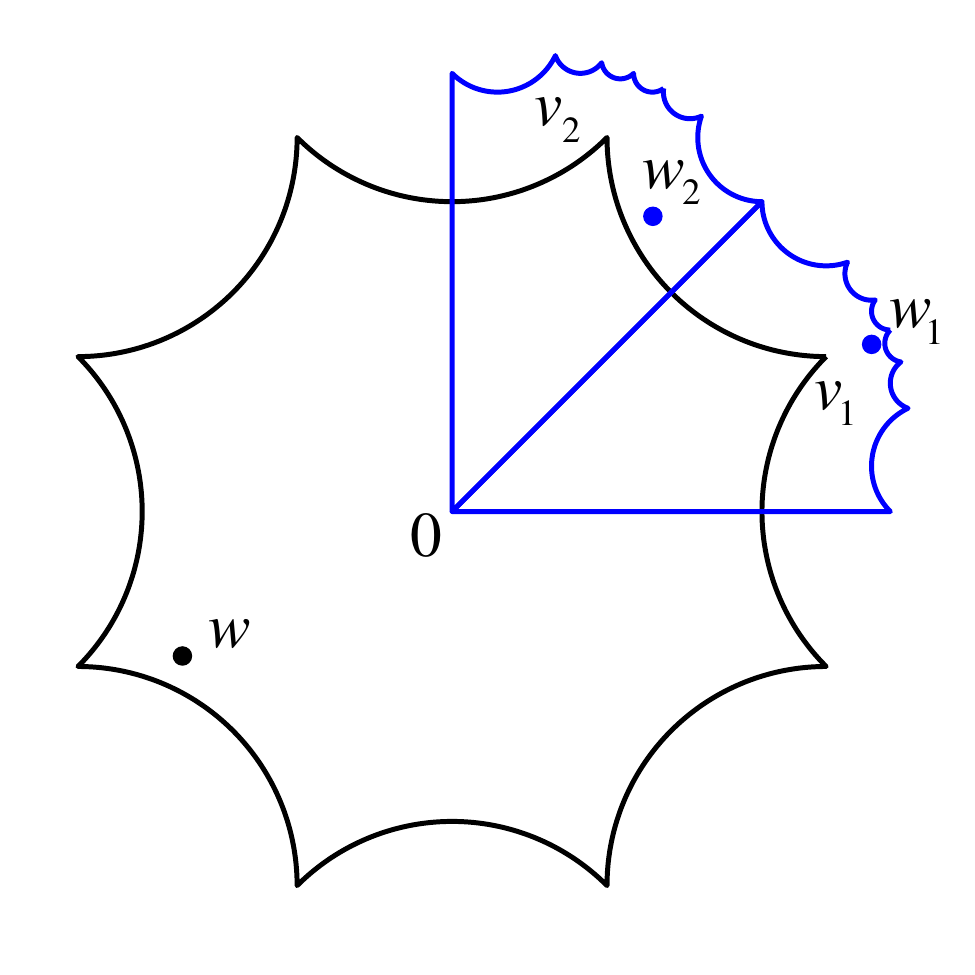}
        \includegraphics[width=0.35\textwidth]{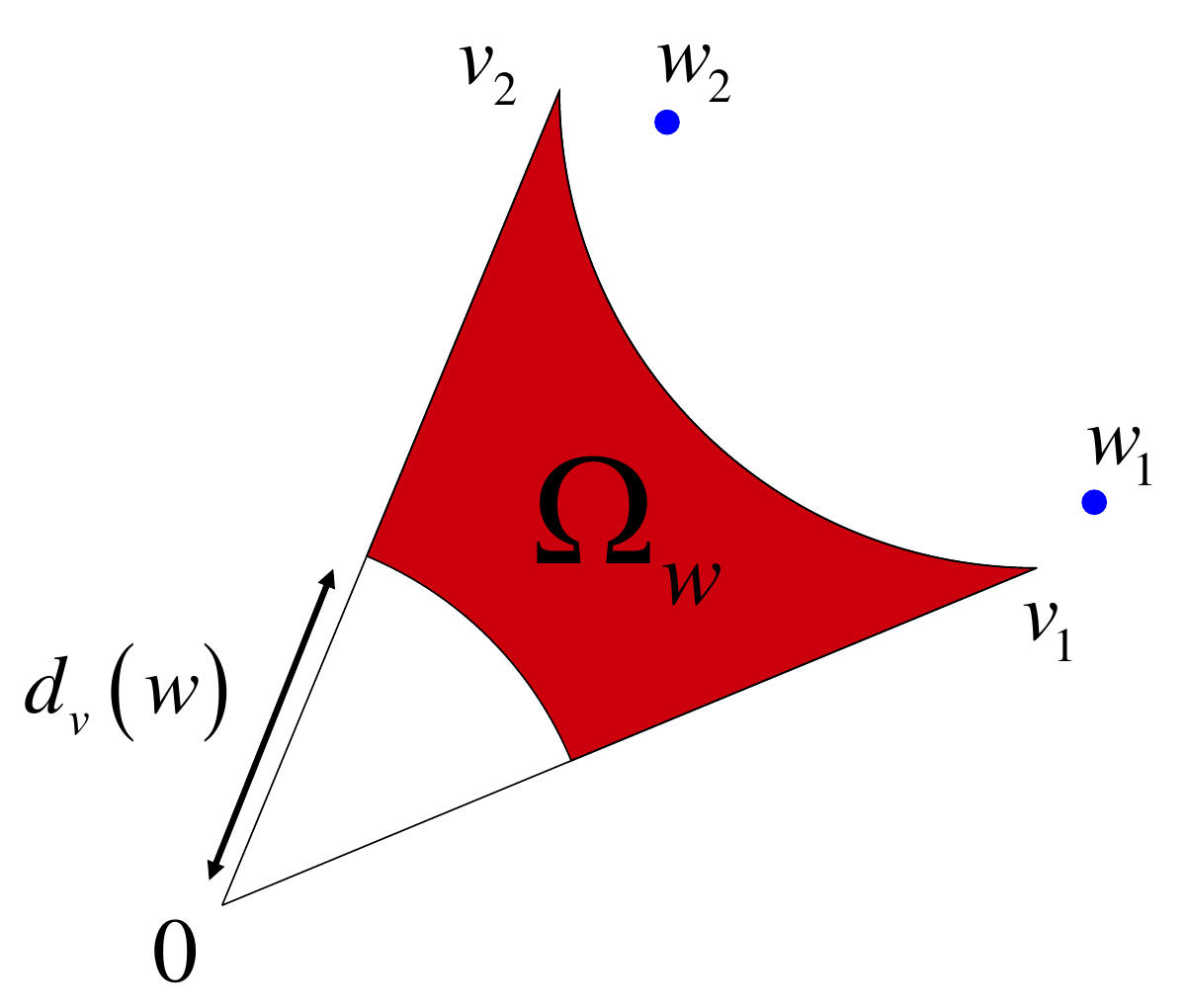}
        \caption{Left: The representatives $w_1$ and $w_2$ of $[w]$ in the dual tessellation 
                 in the ${g=2}$ case. They are images of ${w}$ in the dual polygons 
                 centered at ${v_1}$ and ${v_2}$. 
                 Right: The domain $\Omega_w$ is the set of points ${z}$ 
                 satisfying the conditions in Prop.~\ref{P2} for a given point ${w}$.}
        \label{F3}
    \end{center}
\end{figure}

The rotational and dual symmetries of $S_g$ 
allow us to reduce drastically the set of all the possibilities of where 
the pairs of points $z$ and $w$ can lie to the following subset: 

\begin{proposition} \label{P2}
    To prove Theorem~\ref{T1}, it suffices to consider only the pairs of points 
    $z$ and $w$ satisfying the following conditions: 
    \begin{enumerate}
        \item $z$ lies in triangle ${T = \Delta 0 v_1 v_2}$, Fig.~\ref{F3},
        \item $d_0(z) \geq d_v(w)$,
        \item $d_v(w) \leq s/2$.
    \end{enumerate}
\end{proposition}

Assuming, thanks to Prop.~\ref{P2}, that the first point $z$ lies in the triangle~$T$ 
with vertices $O,v_1,v_2$, let ${w_1,w_2 \in [w]}$ be the two dual-tessellation representatives of $[w]$ 
lying in the two dual polygons centered at $v_1,v_2$, respectively, Fig.~\ref{F3}. 
To prove our desired result ${\delta^\star ([z],[w]) \leq R}$, it suffices to show that 
\begin{equation} \label{Ew12}
    \min \left[ \delta(z,w_1), \delta(z,w_2) \right] \leq R ,
\end{equation}
simply because 
\begin{align}
    \nonumber
    \delta^\star ([z],[w]) = 
    & \min \{ \delta(z^\prime, w^\prime), z^\prime \in [z], w^\prime \in [w] \} \leq \\
    & \min \left[ \delta(z,w_1), \delta(z,w_2) \right] .
\end{align}

Before we prove Eq.~(\ref{Ew12}) using algebra, we take one final 
step to reduce the number of possibilities in locations of $z,w$ to consider. Given $w$ such that 
$d_v(w) \leq s/2$, we define the domain $\Omega_w$ as the set of all $z \in T$ for which $z$ and $w$ 
meet the conditions listed in Prop.~\ref{P2}, Fig.~\ref{F3}. We will prove the following: 

\begin{theorem} \label{T3}
    Given $w$ for which $d_v(w) \leq s/2$, the function 
    \begin{equation}
        \nonumber
        f(z) = \min \left[ \delta(z,w_1), \delta(z,w_2) \right]
    \end{equation}
    attains its global maximum over the domain $\Omega_w$ 
    on its boundary ${\partial \Omega_w}$. 
\end{theorem}

This theorem allows us to reduce the possibilities for the location of $z$ even further, 
to include only the points on the boundary of $\Omega_w$. 

The final step in the proof of Theorem~\ref{T1} is to show that Eq.~(\ref{Ew12}) 
holds for all of the $z,w$ pairs satisfying all the conditions above:

\begin{theorem} \label{T4}
    Given $w$ for which $d_v(w) \leq s/2$, for 
    ${\forall z \in \partial \Omega_w}$ we have 
    \begin{equation}
        \nonumber
        \min \left[ \delta(z,w_1), \delta(z,w_2) \right] \leq R .
    \end{equation}
\end{theorem}

Thus we will have shown that $\mathscr{D}_g \leq R$, and the combination of 
this result with that of Theorem~\ref{T2} proves Theorem~\ref{T1}. 

One important question is how generalizable our proof strategy outlined in this section 
is to other surfaces, such as those appearing in the Poincar\'e theorem~\cite{Katok} or the 
canonical surfaces~\cite{Abrams}. In these surfaces, the gluing identifies not the opposite sides 
of a polygon, but its alternating sides: 1 to 3, 2 to 4, 5 to 7, 6 to 8, and so on. One could 
hope that if the polygon is still a regular ${4g}$-gon, then this gluing would lead to a surface 
with the same diameter ${R}$. However, this is not true---the diameter of these surfaces appears 
to be greater than ${R}$ in simulations. This observation is a reflection of the fact that the 
symmetries of a surface play 
a crucial role in defining its geometric properties including the diameter, so one should select very carefully a right set of 
symmetries of a surface in calculating its diameter. Our proof is an example of this strategy 
applied to a particular class of highly symmetric surfaces. Its generalization to other surfaces 
with other groups of symmetry is an interesting open problem.

The next section contains the complete proofs of all the theorems and propositions above.

\section{Proofs}\label{S4}

\begin{customthm}{2}
    The diameter $\mathscr{D}_g$ of $S_g$ is at least $R$, 
    $\mathscr{D}_g \geq R$. 
\end{customthm}

\noindent \textbf{Proof:} Take ${[z]=[0]}$ and $[w]=[v]$. We will show that 
$\delta^\star([0],[v]) = R$, from which the theorem follows immediately. For 
this, we use that the fundamental polygon $P$ is also the Dirichlet polygon of $[0]$, 
\begin{equation}
    D_0 = \{ p \mid \delta(0,p) < \delta(f(0),p) \text{ for } \forall f \in \mathscr{F}_g, ~ f \neq I \} .
\end{equation}
where ${I}$ is the identity. 
To see this, first note that points $p$ on the geodesic segment 
${\overline{v_iv_{i+1}}}$ satisfy ${\delta(0,p) = \delta(t_i(0),p)}$ owing to the 
symmetry of the fundamental polygon. 
Next, denote by ${H_i}$ the half-plane bounded by this line and 
containing $0$. Then 
\begin{align}
    \nonumber
    & D_0 = 
    \{ p \mid \delta(0,p) < \delta(f(0),p) \text{ for } \forall f \in \mathscr{F}_g, ~ f \neq I \} \subseteq \\
    & \{ \delta(0,p) < \delta(t_i(0),p), ~ 0 \leq i < 4g \} = 
    \cap_{0 \leq i < 4g} H_i ,
\end{align}
and the last expression is just $P$, thus 
${D_0 \subseteq P}$. Since Dirichlet polygons are fundamental polygons, 
which all have the same hyperbolic area (Theorem~9.1.3 in~\cite{Beardon}), we 
have ${\text{area} (D_0) = \text{area} (P)}$ and conclude ${D_0 = P}$. 

Now, since the elements of $[v]$ contained in $P$, which are just the $v_i$, are 
the closest $v$-images to $0$ and lie a distance $R$ from $0$, we have 
$\delta^\star([0],[v]) = R$. Therefore ${\mathscr{D}_g \geq R}$. 

\begin{customprop}{1}
    Let ${\phi}$ be any isometry of $\mathbb{D}$ that maps the fundamental polygon to 
    any of its dual polygons. Then 
    the projection ${\phi^\star : S_g \rightarrow S_g}$ given by 
    \begin{equation}
        \phi^\star ([z]) = [\phi(z)]
    \end{equation}
    is well-defined and is an isometric automorphism of $S_g$. In particular, 
    dual polygons are fundamental domains. 
\end{customprop}

\noindent \textbf{Proof:} 
First we will show the rotational symmetry of the quotient surface. It suffices to show 
that ${\theta \mathscr{F}_g \theta^{-1} = \mathscr{F}_g}$ by the argument given in 
Sec.~\ref{S3}, where ${\theta}$ is the 
counterclockwise rotation of ${\mathbb{D}}$ by ${\pi/2g}$. We will use the 
property of the generators of ${\mathscr{F}_g}$ that ${\theta t_i \theta^{-1} = t_{i+1}}$ for all ${i}$, 
and the fact which follows from elementary group theory that 
${\theta \mathscr{F}_g \theta^{-1} = \mathscr{F}_g}$ if and only if 
${\theta \mathscr{F}_g \theta^{-1} \subseteq \mathscr{F}_g}$. 
We further note that it suffices to show ${\theta t_i \theta^{-1} \subseteq \mathscr{F}_g}$ for 
${\forall i}$, since this implies for 
${\forall f = t_{i_1} t_{i_2} \cdot \cdot \cdot t_{i_k} \in \mathscr{F}_g}$ that 
\begin{equation}
    \theta f \theta^{-1} = 
    \theta t_{i_1} t_{i_2} \cdot \cdot \cdot t_{i_k} \theta^{-1} =
    (\theta t_{i_1} \theta^{-1}) (\theta t_{i_2} \theta^{-1}) \cdot \cdot \cdot (\theta t_{i_k} \theta^{-1}) \in 
    \mathscr{F}_g .
\end{equation}
For an arbitrary generator ${t_i \in \mathscr{F}_g}$ we have
\begin{equation}
    \theta t_i \theta^{-1} = 
    t_{i+1} \in \mathscr{F}_g ,
\end{equation}
and therefore ${\theta \mathscr{F}_g \theta^{-1} = \mathscr{F}_g}$. 
It follows that the isometry ${\theta^\star}$, given by ${\theta^\star ([p]) = [\theta(p)]}$ 
for points ${[p]}$ on the quotient surface, is a well-defined isometry and therefore a symmetry of the 
quotient surface. 

Next, we will show the same for ${\phi^\star}$. Suppose that ${\phi}$ maps the fundamental polygon 
to the dual polygon ${D}$ centered at the vertex ${v}$. Then for any generator ${t_i}$, the map 
${t_i^\prime = \phi t_i \phi^{-1}}$ maps ${D}$ to one of its neighboring dual polygons centered at 
vertex ${v^\prime}$, such that ${t^\prime(v) = v^\prime}$ and ${t^\prime}$ fixes the geodesic line 
between ${v}$ and ${v^\prime}$. 

Now, choose ${h \in \mathscr{F}_g}$ that 
maps ${v}$ to ${v^\prime}$. We will be done if we can show ${t^\prime = h}$, 
since this implies ${t^\prime \in \mathscr{F}_g}$ as desired. 
Since ${h}$ is orientation-preserving, we can write it as a product of 
${t^\prime}$ and an arbitrary rotation by ${\alpha}$ around ${v^\prime}$. 
Then ${h^{-1}}$ is the product of ${(t^\prime)^{-1}}$ and a rotation by 
${- \alpha}$ around ${v}$. 
Pick one of the representatives of the fundamental polygon that has 
${\overline{v v^\prime}}$ as an edge, such that the directed line segment ${\overline{v v^\prime}}$ 
is oriented counterclockwise along the polygon boundary. 
Then by the rotational symmetry of ${\mathscr{F}_g}$, 
every similar element of ${\mathscr{F}_g}$ that maps a vertex 
to its neighbor in the counterclockwise direction has the same form, a translation 
followed by a rotation by ${\alpha}$. Similarly, such maps in 
the clockwise direction are translations followed by rotations by ${- \alpha}$. 

However, we could have chosen the other polygon that has edge ${\overline{v v^\prime}}$ in 
the clockwise orientation, and drawn the same conclusions above with the orientations 
reversed. So we are forced to conclude ${\alpha = - \alpha}$, thus 
either ${\alpha = 0}$ or ${\alpha = \pi}$. But ${h}$ cannot have fixed points, 
and if ${\alpha = \pi}$ then ${h}$ fixes the midpoint between ${v}$ and ${v^\prime}$, 
hence ${\alpha = 0}$, so ${t^\prime = h}$. Therefore ${\phi \mathscr{F} \phi^{-1} = \mathscr{F}}$, 
so ${\phi^\star}$ is a symmetry of the quotient surface. 

\begin{customprop}{2}
    To prove Theorem~\ref{T1}, it suffices to consider only the pairs of points 
    $z$ and $w$ satisfying the following conditions: 
    \begin{enumerate}
        \item $z$ lies in triangle ${T = \Delta 0 v_1 v_2}$, Fig.~\ref{F3},
        \item $d_0(z) \geq d_v(w)$,
        \item $d_v(w) \leq s/2$.
    \end{enumerate}
\end{customprop}

\noindent \textbf{Proof:} For pairs ${z,w}$ that do not meet at least one of the conditions above, 
our strategy will be to find another pair ${z^\prime,w^\prime}$ that do meet the conditions, such 
that ${\delta(z,w) = \delta(z^\prime,w^\prime)}$. Then the diameter of ${S_g}$, which is the maximum of 
${\delta}$ over all ${z,w}$, is equal to the maximum of ${\delta}$ over just the ${z,w}$ which satisfy the 
conditions above, so it will suffice to consider only those pairs. 

\begin{figure}[t!]
    \begin{center}
        \includegraphics[width=0.35\textwidth]{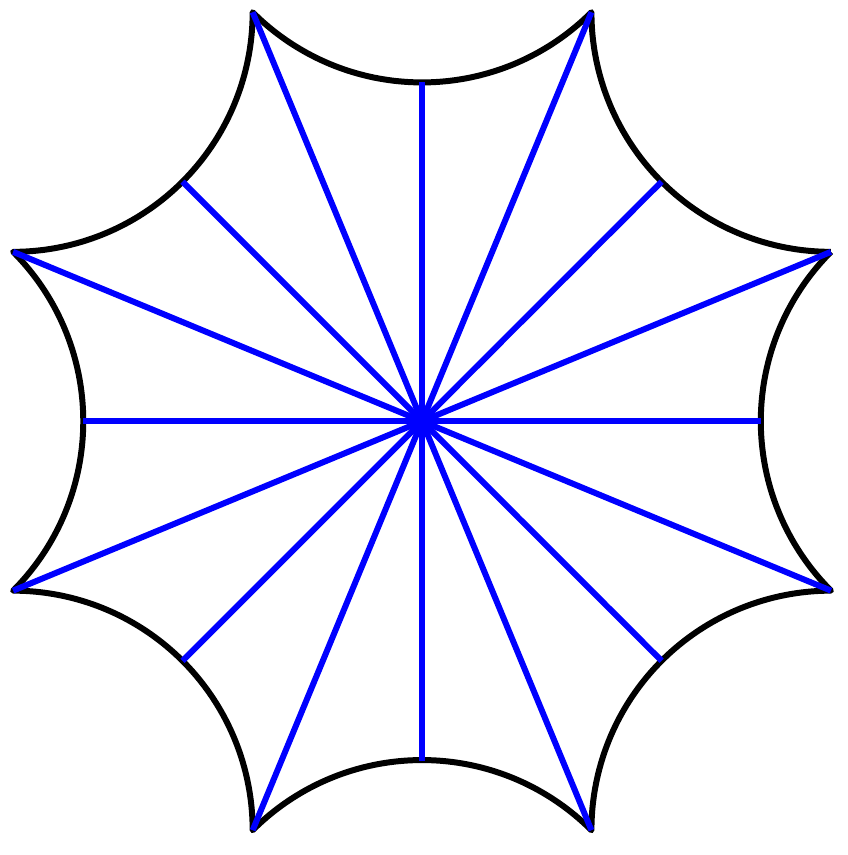}
        \caption{Division of the fundamental polygon into ${8g}$ isosceles triangles for ${g=2}$.}
        \label{F4}
    \end{center}
\end{figure}

Starting with an arbitrary pair ${z,w}$, 
let $\phi$ be any of the isometries from the fundamental 
polygon to a dual polygon, and $\phi^\star$ the corresponding isometry of the quotient surface. 
Note that $\phi^\star$ takes ${[0] \rightarrow [v]}$ and ${[v] \rightarrow [0]}$. 
Therefore ${d_0(\phi(p)) = d_v(p)}$ and ${d_v(\phi(p)) = d_0(p)}$ for points ${p \in \mathbb{D}}$. 
Now, for the given $z,w$ consider the four quantities ${d_0(z), d_v(z), d_0(w), d_v(w)}$. 
We can divide the fundamental polygon into ${8g}$ congruent isosceles triangles (Fig.~\ref{F4}) each with 
one vertex at ${0}$, one vertex at some ${v_i}$, and one vertex at the midpoint of an edge of the polygon, 
and side lengths ${R,s/2,s/2}$. Since ${z}$ lies in one of these isosceles triangles, 
we have ${d_0(z) + d_v(z) \leq s/2 + s/2 = s}$, and similarly for ${w}$. It follows that 
\begin{equation}
    M = \min \{ d_0(z), d_v(z), d_0(w), d_v(w) \} \leq s/2 .
\end{equation}
We now consider four cases, and in each case construct ${z^\prime,w^\prime}$ for which 
conditions (2) and (3) are satisfied and ${\delta(z^\prime,w^\prime) = \delta(z,w)}$. 
\begin{enumerate}
    \item ${M = d_0(z)}$. Let ${z^\prime = \phi(w)}$ and ${w^\prime = \phi(z)}$. Then 
        \begin{align}
            d_0 (z^\prime) = d_v(w) & \geq d_0(z) = d_v(w^\prime) \text{ and} \\
            d_v(w^\prime) = d_0(z) & \leq s/2 .
        \end{align}
    \item ${M = d_v(z)}$. Let ${z^\prime = w}$ and ${w^\prime = z}$. Then 
        \begin{align}
            d_0 (z^\prime) = d_0(w) & \geq d_v(z) = d_v(w^\prime) \text{ and} \\
            d_v(w^\prime) = d_v(z) & \leq s/2 .
        \end{align}
    \item ${M = d_0(w)}$. Let ${z^\prime = \phi(z)}$ and ${w^\prime = \phi(w)}$. Then 
        \begin{align}
            d_0 (z^\prime) = d_v(z) & \geq d_0(w) = d_v(w^\prime) \text{ and} \\
            d_v(w^\prime) = d_0(w) & \leq s/2 .
        \end{align}
    \item ${M = d_v(w)}$. Let ${z^\prime = z}$ and ${w^\prime = w}$. Then 
        \begin{align}
            d_0 (z^\prime) = d_0(z) & \geq d_v(w) = d_v(w^\prime) \text{ and} \\
            d_v(w^\prime) = d_v(w) & \leq s/2 .
        \end{align}
\end{enumerate}

We have already shown the rotational symmetry of ${S_g}$ in 
the proof of Prop.~\ref{P1}. Starting with $z^\prime,w^\prime$ and rotating a sufficient 
number of times by ${\pi/2g}$, we can place ${z^\prime}$ in $T$. 
Furthermore, since rotations map vertices to vertices and fix ${0}$, 
they preserve the quantities ${d_0(z^\prime), d_v(z^\prime), d_0(w^\prime), d_v(w^\prime)}$, 
so the rotated ${z^\prime,w^\prime}$ now satisfy all three conditions. 

\begin{customthm}{3}
    Given $w$ for which $d_v(w) \leq s/2$, the function 
    \begin{equation}
        \nonumber
        f(z) = \min \left[ \delta(z,w_1), \delta(z,w_2) \right]
    \end{equation}
    \noindent attains its global maximum over the domain $\Omega_w$ 
    on its boundary ${\partial \Omega_w}$. 
\end{customthm}

\noindent \textbf{Proof:} Since ${f(z)}$ is continuous and bounded on the 
compact region ${\Omega_w \cup \partial \Omega_w}$, 
$f(z)$ attains a global maximum on this region. 

Let $M$ denote the midline of $w_1$ and $w_2$ (the locus of points equidistant 
from $w_1$ and $w_2$, also the perpendicular bisector of the geodesic segment $\overline{w_1w_2}$). 
Then $M$ divides $\Omega_w$ into some set of open regions, $\{S_i\}$. 
We will first show that for any ${i}$ and ${z_0 \in S_i}$, 
$z_0$ cannot be a global maximum of $f(z)$. 

Inside ${S_i}$ we have ${\delta(z, w_1) \neq \delta(z, w_2)}$, and 
therefore either ${\delta(z, w_1) < \delta(z, w_2)}$ or ${\delta(z, w_2) < \delta(z, w_1)}$. 
Assume without loss of generality 
that ${\delta(z, w_1) < \delta(z, w_2)}$. Then ${f(z_0) = \delta(z_0, w_1)}$ and 
${f(z) = \delta(z, w_1)}$ in some open set around $z_0$. Therefore, 
we can always increase the value of $f$ by moving $z_0$ in some suitable direction 
(away from $w_1$). 

It remains to consider ${z_0 \in M \cap \Omega_w}$. The set~${M \cap \Omega_w}$ is a union of segments 
${\{M_i\}}$ which lie 
in $\Omega_w$ and have endpoints on $\partial \Omega_w$. Take ${z_0 \in M_i}$. 
In this case we also have ${f(z_0) = \delta(z_0, w_1)}$. It suffices to show 
that $\delta(z, w_1)$ is maximized at one of the endpoints of $M_i$. However, this is 
a well-known result on the distance from points to line segments 
(see Theorem 2.2 in~\cite{Ramsay}). Since these endpoints are on $\partial \Omega_w$, 
we have shown that the global maximum of ${f(z)}$ is attained on $\partial \Omega_w$. 

\begin{customthm}{4}
    Given $w$ for which $d_v(w) \leq s/2$, for 
    ${\forall z \in \partial \Omega_w}$ we have 
    \begin{equation}
        \nonumber
        \min \left[ \delta(z,w_1), \delta(z,w_2) \right] \leq R .
    \end{equation}
\end{customthm}

\noindent \textbf{Proof:} We abbreviate ${a = d_0(z)}$, ${b = d_v(w)}$, 
${R^{\prime} = \tanh \left( \frac{R}{2} \right)}$, and ${s^{\prime} = \tanh \left( \frac{s}{4} \right)}$. 

\noindent \\\textbf{CASE I:} We look at the radial segments of ${\partial \Omega_w}$ first, see 
Fig.~\ref{F3}, right. For ${z}$ on ${\overline{0v_0}}$, 
consider the triangle with vertices ${z}$, ${v_0}$, and ${w_1}$ and side lengths ${b}$, ${\delta(z,w_1)}$, and 
${R - a}$, respectively. By the triangle inequality, 
\begin{align}
    \nonumber
    & b + (R - a) \geq \delta(z,w_1) \Rightarrow \\
    & \delta(z,w_1) \leq R - (a - b) \leq R,
\end{align}
and a similar argument for ${z}$ on ${\overline{0v_2}}$ gives ${\delta(z,w_2) \leq R}$, so that 
in either case ${\min \left[ \delta(z,w_1), \delta(z,w_2) \right] \leq R}$. 

\begin{figure}[t!]
    \begin{center}
        \includegraphics[width=0.35\textwidth]{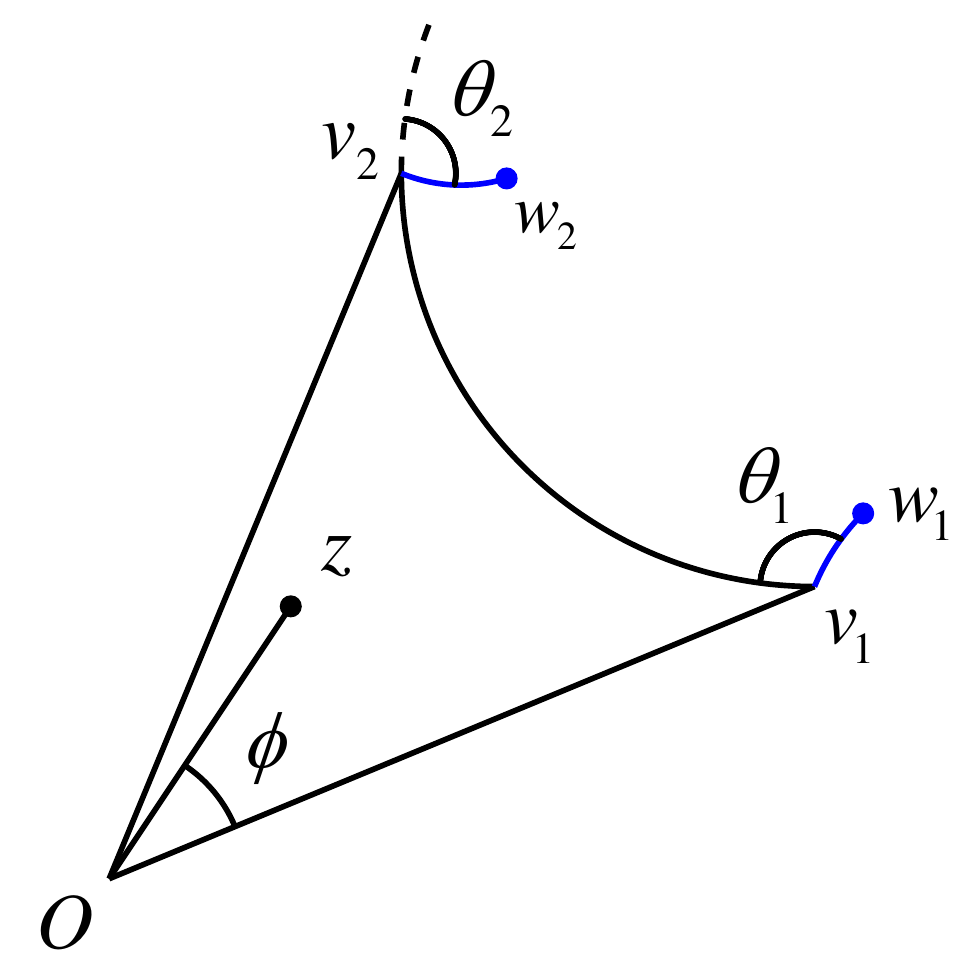}
        \caption{Illustration of ${z}$, ${w_1}$, and ${w_2}$, and angles 
                 ${\phi}$, ${\theta_1}$, ${\theta_2}$.}
        \label{F5}
    \end{center}
\end{figure}

\noindent \\\textbf{CASE II:} Next we consider the circular arc of radius ${b}$ centered at the origin 
(Fig.~\ref{F3}, right). For ${z}$ on this arc we have ${a = b}$. Let ${\phi}$ denote the 
angle between ${\overline{0z}}$ and ${\overline{0v_1}}$, so that ${\phi \in \left[ 0, \frac{\pi}{2g} \right]}$. 
Let ${\theta_1}$ denote the angle between 
${\overline{v_1w_1}}$ and ${\overline{v_1v_2}}$, and ${\theta_2}$ denote the angle between 
${\overline{v_2w_2}}$ and 
the extension of ${\overline{v_1v_2}}$, Fig.~\ref{F5}. 

We have the constraint ${\theta_1 = \theta_2}$, i.e. $w_1$ and $w_2$ make the 
same angle with the axis ${\overline{v_1v_0}}$, which is a consequence of Prop.~\ref{P1}. 
Now, we assume that ${\delta(z,w_1) > R}$ and 
${\delta(z,w_2) > R}$ and seek a contradiction. 
Using the hyperbolic sine and cosine theorems to calculate 
the distances, we get after some algebra
\begin{align} \label{E26}
    & \nonumber
    \delta(z,w_1) > R \Leftrightarrow \\
    & \nonumber
    \left( 1 + R^{\prime ^2} \right) 
    - 2 R^{\prime} \coth a \left[ \cos \left( \theta_1 + \frac{\pi}{4g} \right) 
    + \cos \phi \right] + \\
    &
    \cos \left( \theta_1 + \phi + \frac{\pi}{4g} \right)
    + R^{\prime ^2} \cos \left( \theta_1 - \phi + \frac{\pi}{4g} \right)
    > 0,
\end{align}
and a similar equation arises from ${\delta(z,w_2) > R}$. 
Solving Eq.~(\ref{E26}) yields 
\begin{equation}
    \frac{\cos \left(\frac{\theta_1 + \phi}{2} + \frac{\pi}{8g} \right)}{\sin \left( \frac{\theta_1 - \phi}{2} + \frac{\pi}{8g} \right)} \in
    \left( - \infty, R^{\prime} \tanh \left( \frac{a}{2} \right) \right) \cup
    \left( \frac{R^{\prime}}{\tanh \left( \frac{a}{2} \right)}, \infty \right) .
\end{equation}
Notice that ${R^{\prime} \tanh \left( \frac{a}{2} \right)}$ is an 
increasing function of ${a}$, while ${\frac{R^{\prime}}{\tanh \left( \frac{a}{2} \right)}}$ 
is a decreasing function of ${a}$, so it suffices to consider the case where 
${a}$ is maximal, ${a = \frac{s}{2}}$. Plugging this in gives 
\begin{equation}
    \frac{\cos \left(\frac{\theta_1 + \phi}{2} + \frac{\pi}{8g} \right)}{\sin \left( \frac{\theta_1 - \phi}{2} + \frac{\pi}{8g} \right)} \in
    \left( - \infty, R^{\prime} s^\prime \right) \cup
    \left( \frac{R^{\prime}}{s^\prime}, \infty \right) .
\end{equation}
Solving for ${\alpha = \theta_1 - \phi + \frac{\pi}{4g}}$ gives
\begin{equation} \label{E29}
    \tan \left( \frac{\alpha}{2} \right) \in
    \left( \frac{\cos \phi - R^\prime/s^\prime}{\sin \phi},
    \frac{\cos \phi - R^\prime s^\prime}{\sin \phi} \right) .
\end{equation}
A similar analysis of the second inequality, ${\delta(z,w_2) > R}$ using 
${\theta_1 = \theta_2}$ gives 
\begin{align} \label{E30}
    \nonumber
    \tan \left( \frac{\alpha}{2} \right) \in
    & \left( - \infty, 
    \frac{\sin \left( \frac{\pi}{2g} - \phi \right)}{\cos \left( \frac{\pi}{2g} - \phi \right) - R^{\prime}/s^{\prime}}
    \right) \cup \\ & \left(
    \frac{\sin \left( \frac{\pi}{2g} - \phi \right)}{\cos \left( \frac{\pi}{2g} - \phi \right) - R^{\prime} s^{\prime}},
    \infty \right) .
\end{align}
However, since the inequalities
\begin{equation}
    \frac{\cos \phi - R^\prime/s^\prime}{\sin \phi} \leq
    \frac{\sin \left( \frac{\pi}{2g} - \phi \right)}{\cos \left( \frac{\pi}{2g} - \phi \right) - R^{\prime}/s^{\prime}}
\end{equation}
and
\begin{equation}
    \frac{\cos \phi - R^\prime s^\prime}{\sin \phi} \geq
    \frac{\sin \left( \frac{\pi}{2g} - \phi \right)}{\cos \left( \frac{\pi}{2g} - \phi \right) - R^{\prime} s^{\prime}}
\end{equation}
have no solutions in ${\phi \in \left[ 0, \frac{\pi}{2g} \right]}$, 
the sets in Eqs.~(\ref{E29},\ref{E30}) are disjoint, and we conclude that the system 
\begin{equation}
    \begin{cases}
        \delta(z,w_1) > R \\
        \delta(z,w_2) > R
    \end{cases}
\end{equation}
has no solutions, so ${\min \left[ \delta(z,w_1), \delta(z,w_2) \right] \leq R}$. 

\noindent \\\textbf{CASE III:} Finally, we consider ${z}$ on the edge ${\overline{v_1v_2}}$. 
Let ${x = \delta(z,v_1)}$ and ${\alpha = \angle v_2v_1w_1}$. 
Applying the first hyperbolic law of cosines to triangle ${\Delta z v_1 w_1}$ gives
\begin{equation}
    \cosh \left( \delta(z, w_1) \right) = \cosh b \cosh x - \sinh b \sinh x \cos \alpha .
\end{equation}
Assuming for the sake of contradiction that ${\delta(z,w_1) > R}$ and 
${\delta(z,w_2) > R}$, then from the first inequality we find
\begin{equation}
    \cosh b \cosh x - \sinh b \sinh x \cos \alpha > \cosh R ,
\end{equation}
which implies
\begin{equation} \label{E36}
    \cos \alpha < \frac{\cosh b \cosh x - \cosh R}{\sinh b \sinh x} .
\end{equation}
Similarly, the second inequality ${\delta(z,w_2) > R}$ yields
\begin{equation} \label{E37}
    \cos \alpha > \frac{\cosh R - \cosh b \cosh (s - x)}{\sinh b \sinh (s - x)} .
\end{equation}
However, for inequalities~(\ref{E36}) and~(\ref{E37}) to have solutions in ${\alpha}$, 
we must have
\begin{equation}
    \frac{\cosh b \cosh x - \cosh R}{\sinh b \sinh x} >
    \frac{\cosh R - \cosh b \cosh (s - x)}{\sinh b \sinh (s - x)} .
\end{equation}
Simplifying yields
\begin{equation}
    \cosh \left( \frac{s}{2} \right) \cosh \left( x - \frac{s}{2} \right) <
    \cosh b .
\end{equation}
Notice that the right hand side is an increasing function of ${b}$, thus 
it suffices to disprove the inequality 
in the case where ${b}$ is maximal, ${b = \frac{s}{2}}$. Plugging this in and simplifying 
the resulting inequality gives
\begin{equation}
    \cosh \left( x - \frac{s}{2} \right) < 1 ,
\end{equation}
which has no solutions in ${x}$. Therefore our assumption was false, so 
${\min \left[ \delta(z,w_1), \delta(z,w_2) \right] \leq R}$.

\begin{customthm}{1}
    ${\mathscr{D}_g = R}$. 
\end{customthm}

\noindent \textbf{Proof:} We first show that for ${\forall [z],[w]}$, 
\begin{equation} \label{E51}
    \delta^\star([z],[w]) \leq R .
\end{equation}
By Prop.~\ref{P2}, it suffices to consider the case ${z \in \Delta 0v_1v_2}$, 
${d_0(z) \geq d_v(w)}$, ${d_v(w) \leq s/2}$. Then by Theorem~\ref{T3}, for all 
such $w$ the function ${f(z) = \delta^\star([z],[w])}$ defined on the domain 
${\Omega_w}$ attains a global maximum $M$ on the boundary ${\partial \Omega_w}$, and 
by Theorem~\ref{T4}, $M \leq R$. This proves Eq.~(\ref{E51}), 
and therefore $\mathscr{D}_g \leq R$. But Theorem~\ref{T2} says that $\mathscr{D}_g \geq R$.  
We conclude that $\mathscr{D}_g = R$. 

\noindent \\\\\textbf{ACKNOWLEDGEMENTS}

\noindent \\We thank Florent Balacheff, Maxime Fortier Bourque, Antonio Costa, 
Benson Farb, Svetlana Katok, Gabor Lippner, Marissa Loving, Curtis McMullen, 
Hugo Parlier, John Ratcliffe, and Chaitanya Tappu for 
useful discussions and suggestions. This work was supported by NSF grant 
Nos.~IIS-1741355 and~CCF-2311160.

\bibliography{sn-bibliography}
\nocite{*}

\end{document}